\magnification=1200
%%%%%%%%%%%%%%%%%% Macro Definitions %%%%%%%%%%%%%%%%%

\def\title#1{{\titlefont\noindent #1\bigskip}}

\def\author#1{{\largefont\noindent #1}\medskip}

\def\beginlinemode{\endmode
 \begingroup\obeylines\def\endmode{\par\endgroup}}
\let\endmode=\par

\newbox\theaddress
\def\address{\smallskip\beginlinemode\parindent 0in\getaddress}
{\obeylines
\gdef\getaddress #1 
 #2
 {#1\gdef\addressee{#2}%
   \global\setbox\theaddress=\vbox\bgroup\raggedright%
    \everypar{\hangindent2em}#2
   \def\endaddress{\egroup\endgroup \copy\theaddress \medskip}}}

\def\thanks#1{\footnote{}{\eightpoint #1}}

\long\def\Abstract#1{{\eightpoint\narrower\vskip\baselineskip\noindent
#1\smallskip}}

\def\skipfirstword#1 {}

\def\ir#1{\csname #1\endcsname}

\newdimen\currentht
\newbox\droppedletter
\newdimen\droppedletterwdth
\newdimen\drophtinpts
\newdimen\dropindent

\def\irrnSection#1#2{\edef\tttempcs{\ir{#2}}
\vskip\baselineskip\penalty-3000
{\largefont\bf\noindent \expandafter\skipfirstword\tttempcs. #1}
\vskip6pt}

\def\irSubsection#1#2{\edef\tttempcs{\ir{#2}}
\vskip\baselineskip\penalty-3000
{\bf\noindent \expandafter\skipfirstword\tttempcs. #1}
\vskip6pt}

\def\irSubsubsection#1#2{\edef\tttempcs{\ir{#2}}
\vskip\baselineskip\penalty-3000
{\noindent \expandafter\skipfirstword\tttempcs. #1}
\vskip6pt}

\def\References{\vbox to.25in{\vfil}\noindent{}{\bf References}
\vskip6pt\par}

\def\References{\vskip6pt\noindent{}{\bf References}
\vskip6pt\par}

\def\baselinebreak{\par \ifdim\lastskip<6pt
         \removelastskip\penalty-200\vskip6pt\fi}

\long\def\prclm#1#2#3{\baselinebreak
\noindent{\bf \csname #2\endcsname}:\enspace{\sl #3\par}\baselinebreak}

\def\rem#1#2{\baselinebreak\noindent{\bf \csname #2\endcsname}:}

\def\qed{{$\diamondsuit$}\vskip6pt}

\def\bibitem#1{\par\indent\llap{\rlap{\bf [#1]}\indent}\indent\hangindent
2\parindent\ignorespaces}

\long\def\eatit#1{}

\def\leftheadlinetext{}
\def\rightheadlinetext{}

\def\leftheadline{{\eightrm\folio\hfil \leftheadlinetext\hfil}}
\def\rightheadline{{\eightrm\hfil\rightheadlinetext\hfil\folio}}

\headline={\ifnum\pageno=1\hfil\else
\ifodd\pageno\rightheadline\else\leftheadline\fi\fi}

\def\tenpoint{\def\rm{\fam0\tenrm}
\textfont0=\tenrm \scriptfont0=\sevenrm \scriptscriptfont0=\fiverm
\textfont1=\teni \scriptfont1=\seveni \scriptscriptfont1=\fivei
\def\mit{\fam1} \def\oldstyle{\fam1\teni}
\textfont2=\tensy \scriptfont2=\sevensy \scriptscriptfont2=\fivesy
\def\cal{\fam2}
\textfont3=\tenex \scriptfont3=\tenex \scriptscriptfont3=\tenex
\def\it{\fam\itfam\tenit} % \it is family 4
\textfont\itfam=\tenit
\def\sl{\fam\slfam\tensl} % \sl is family 5
\textfont\slfam=\tensl
\def\bf{\fam\bffam\tenbf} % \bf is family 6
\textfont\bffam=\tenbf \scriptfont\bffam=\sevenbf
\scriptscriptfont\bffam=\fivebf
\def\tt{\fam\ttfam\tentt} % \tt is family 7
\textfont\ttfam=\tentt
\normalbaselineskip=12pt
\setbox\strutbox=\hbox{\vrule height8.5pt depth3.5pt  width0pt}%
\normalbaselines\rm}

\def\eightpoint{\def\rm{\fam0\eightrm}%
\textfont0=\eightrm \scriptfont0=\sixrm \scriptscriptfont0=\fiverm
\textfont1=\eighti \scriptfont1=\sixi \scriptscriptfont1=\fivei
\def\mit{\fam1} \def\oldstyle{\fam1\eighti}%
\textfont2=\eightsy \scriptfont2=\sixsy \scriptscriptfont2=\fivesy
\def\cal{\fam2}%
\textfont3=\tenex \scriptfont3=\tenex \scriptscriptfont3=\tenex
\def\it{\fam\itfam\eightit} % \it is family 4
\textfont\itfam=\eightit
\def\sl{\fam\slfam\eightsl} % \sl is family 5
\textfont\slfam=\eightsl
\def\bf{\fam\bffam\eightbf} % \bf is family 6
\textfont\bffam=\eightbf \scriptfont\bffam=\sixbf
\scriptscriptfont\bffam=\fivebf
\def\tt{\fam\ttfam\eighttt} % \tt is family 7
\textfont\ttfam=\eighttt
\normalbaselineskip=9pt%
\setbox\strutbox=\hbox{\vrule height7pt depth2pt  width0pt}%
\normalbaselines\rm}

\def\largefont{\def\rm{\fam0\largerm}
\textfont0=\largerm \scriptfont0=\largescriptrm \scriptscriptfont0=\tenrm
\textfont1=\largei \scriptfont1=\largescripti \scriptscriptfont1=\teni
\def\mit{\fam1} \def\oldstyle{\fam1\teni}
\textfont2=\largesy %\scriptfont2=\sevensy \scriptscriptfont2=\fivesy
\def\cal{\fam2}
%\textfont3=\largeex %\scriptfont3=\tenex \scriptscriptfont3=\tenex
\def\it{\fam\itfam\largeit} % \it is family 4
\textfont\itfam=\largeit
\def\sl{\fam\slfam\largesl} % \sl is family 5
\textfont\slfam=\largesl
\def\bf{\fam\bffam\largebf} % \bf is family 6
\textfont\bffam=\largebf %\scriptfont\bffam=\sevenbf 
\scriptscriptfont\bffam=\fivebf
\def\tt{\fam\ttfam\largett} % \tt is family 7
\textfont\ttfam=\largett
\normalbaselineskip=17.28pt
\setbox\strutbox=\hbox{\vrule height12.25pt depth5pt  width0pt}%
\normalbaselines\rm}

\def\titlefont{\def\rm{\fam0\titlerm}
\textfont0=\titlerm \scriptfont0=\largescriptrm \scriptscriptfont0=\tenrm
\textfont1=\titlei \scriptfont1=\largescripti \scriptscriptfont1=\teni
\def\mit{\fam1} \def\oldstyle{\fam1\teni}
\textfont2=\titlesy %\scriptfont2=\sevensy \scriptscriptfont2=\fivesy
\def\cal{\fam2}
%\textfont3=\largeex %\scriptfont3=\tenex \scriptscriptfont3=\tenex
\def\it{\fam\itfam\titleit} % \it is family 4
\textfont\itfam=\titleit
\def\sl{\fam\slfam\titlesl} % \sl is family 5
\textfont\slfam=\titlesl
\def\bf{\fam\bffam\titlebf} % \bf is family 6
\textfont\bffam=\titlebf %\scriptfont\bffam=\sevenbf 
\scriptscriptfont\bffam=\fivebf
\def\tt{\fam\ttfam\titlett} % \tt is family 7
\textfont\ttfam=\titlett
\normalbaselineskip=24.8832pt
\setbox\strutbox=\hbox{\vrule height12.25pt depth5pt  width0pt}%
\normalbaselines\rm}

\nopagenumbers

\font\eightrm=cmr8
\font\eighti=cmmi8
\font\eightsy=cmsy8
\font\eightbf=cmbx8
\font\eighttt=cmtt8
\font\eightit=cmti8
\font\eightsl=cmsl8
\font\sixrm=cmr6
\font\sixi=cmmi6
\font\sixsy=cmsy6
\font\sixbf=cmbx6

\font\largerm=cmr12 at 17.28pt
\font\largei=cmmi12 at 17.28pt
\font\largescriptrm=cmr12 at 14.4pt
\font\largescripti=cmmi12 at 14.4pt
\font\largesy=cmsy10 at 17.28pt
\font\largebf=cmbx12 at 17.28pt
\font\largett=cmtt12 at 17.28pt
\font\largeit=cmti12 at 17.28pt
\font\largesl=cmsl12 at 17.28pt

\font\titlerm=cmr12 at 24.8832pt
\font\titlei=cmmi12 at 24.8832pt
\font\titlesy=cmsy10 at 24.8832pt
\font\titlebf=cmbx12 at 24.8832pt
\font\titlett=cmtt12 at 24.8832pt
\font\titleit=cmti12 at 24.8832pt
\font\titlesl=cmsl12 at 24.8832pt

\tenpoint

%%%%%%%%%%%%%% Internal References %%%%%%%%%%%%%%%%%
%Internal references 

%%%%%%%%%%%%%%%%%%%%%%%%%%%%%%% Begin Paper %%%%%%%%%%%%%%%%%%%%%%%%%%%

\def\manyby{\hbox to.75in{\hrulefill}}
\hsize 5.41667in 
\vsize 7.5in

\def\manyby{\hbox to.75in{\hrulefill}}

\tolerance 3000
\hbadness 3000

\def\item#1{\par\indent\indent\llap{\rlap{#1}\indent}\hangindent
2\parindent\ignorespaces}

\def\itemitem#1{\par\indent\indent
\indent\llap{\rlap{#1}\indent}\hangindent
3\parindent\ignorespaces}

\def\refBal{Bal}
\def\refBau{Bau}
\def\refBi{Bi}
\def\refCat{Ca}
\def\refCM{CM}
\def\refDGM{DGM}
\def\refdem{De}
\def\refEL{EL}
\def\refE{Ev}
\def\refGib{Gi}
\def\refHHF{HHF}
\def\refars{H1}
\def\refnagprob{H2}
\def\refsurv{H3}
\def\refHR{HR}
\def\refHi{Hi}
\def\refKu{Ku1}
\def\refKub{Ku2}
\def\refNone{N1}
\def\refNtwo{N2}
\def\reframanujam{Ra}
\def\refRoe{R1}
\def\refRoeb{R2}
\def\refS{S}
\def\refST{ST}
\def\refterakawa{T}
\def\refXa{Xu1}
\def\refXb{Xu2}

\def\C#1{\hbox{$\cal #1$}}

\def\pr#1{\hbox{{\bf P}${}^{#1}$}}
\def\cite#1{[\ir{#1}]}

\def\leftheadlinetext{Brian Harbourne}
\def\rightheadlinetext{Seshadri constants on algebraic surfaces}

\title{Seshadri constants and very ample divisors on algebraic surfaces}

\author{Brian Harbourne}

\address
Department of Mathematics and Statistics
University of Nebraska-Lincoln
Lincoln, NE 68588-0323
email: bharbour@math.unl.edu
WEB: {\tt http://www.math.unl.edu/$\sim$bharbour/}
\smallskip
April 3, 2001\endaddress
\vskip-\baselineskip

\thanks{\vskip -6pt
\noindent This work benefitted from a National Science Foundation grant.
Also, I thank Prof. Szemberg for his comments
and for sending me his recent preprints.
\smallskip
\noindent 2000 {\it Mathematics Subject Classification. } 
%14C20 Divisors, linear systems, invertible sheaves
%14J60 Vector bundles on surfaces and 
%      higher dimensional varieties, and their moduli
%14J26 Rational and ruled surfaces
%14H50 Plane and space curves
%14N05 Projective techniques
Primary 14C20, 14J60, 14N05. 
Secondary 14J26, 14H50.
\smallskip
\noindent {\it Key words and phrases. }  Seshadri constant,
algebraic surface, line bundles, linear systems, 
rational surface, fat points, nef, very ample divisor.\smallskip}

\vskip\baselineskip
\Abstract{Abstract: A broadly applicable geometric approach for 
constructing nef divisors on blow ups of algebraic surfaces at $n$ 
general points is given; it works for all surfaces
in all characteristics for any $n$. This construction is used to obtain
substantial improvements for currently known lower bounds for 
$n$ point Seshadri constants. Remarks are included about a range 
of applications to classical problems involving linear systems on \pr2.}
\vskip\baselineskip

\irrnSection{Introduction}{intro}
This paper presents a broadly applicable geometric approach to
building nef divisors on surfaces. Our main application is to
obtaining bounds on multipoint Seshadri constants
for $n$ general points on surfaces $X$. What we find is that,
for $n$ sufficiently large, all of the main results for
$X=\pr2$ hold for surfaces generally. 

To begin, let $X$ be an algebraic surface (by which we will always mean
a reduced, irreducible, normal projective variety of dimension 2,
over an algebraically closed field of arbitrary characteristic).
Let $L$ be a nef divisor on $X$, let 
$l=L^2$, and let $p_1,\ldots,p_n$ be distinct 
points of $X$. Seshadri constants were introduced in \cite{refdem};
more generally, multiple point versions have 
been studied in \cite{refBau}, \cite{refBi}, \cite{refKu}, \cite{refXb},
\cite{refS} and \cite{refST}. To recall, the multiple point 
Seshadri constant
$\epsilon(L, p_1,\ldots,p_n)$ is defined 
to be the supremum of all rational numbers $\varepsilon$ such that
$\pi^*L-\varepsilon(E_1+\cdots+E_n)$ is a nef ${\bf Q}$-divisor, where 
$\pi:Y\to X$ is the morphism blowing up the points $p_i$, $1\le i\le n$,
and $E_i$ is the exceptional divisor corresponding to $p_i$.
We will often be concerned with finding lower bounds for
$\epsilon(L, p_1,\ldots,p_n)$ which hold on an open set in $X^n$.
Thus, given a lower bound $c$, it will be convenient 
to write $\epsilon(L, n)\ge c$
to mean that $\epsilon(L, p_1,\ldots,p_n)\ge c$
holds on an open set of $n$-tuples of points $p_i$ of $X$.

It is not hard to see that $\epsilon(L, p_1,\ldots,p_n)\le \sqrt{l/n}$
(see \ir{scva}), and, as remarked in \cite{refXb}, it follows over 
${\bf C}$ from \cite{refEL} that $\epsilon(L, p_1,\ldots,p_n)\ge 1$ 
for sufficiently general points $p_i$, if $l>n$.
Another lower bound for $\epsilon(L, p_1,\ldots,p_n)$ for
sufficiently general points $p_i$ over ${\bf C}$  
follows from the main result of \cite{refKu},
but (in dimension 2) this lower bound 
is never more than $\sqrt{l/n}\sqrt{1-1/n}$. 
Our results are of interest mainly when $l<n$: for any given
very ample divisor $L$, our main result, \ir{introThm}, 
obtains (in view of \ir{introProp})
better bounds than $\sqrt{l/n}\sqrt{1-1/n}$ 
for almost all $n$ sufficiently large. 

There are very few cases for which the value of
$\epsilon(L, p_1,\ldots,p_n)$ is known for $n$ general or generic
points $p_i$. One important case, which has seen a great deal
of attention beginning with Nagata's work on Hilbert's 14th Problem, 
is when $L$ is very ample and $l=1$,
which forces $X$ to be \pr2, but even here, $\epsilon(L, n)$
is known only when $n<9$ (see \ir{smalln}) or when
$n$ is a square. In fact, Nagata's conjecture \cite{refNone},
that (in different terminology)
$\epsilon(L, p_1,\ldots,p_n)=1/\sqrt{n}$ should hold 
for $n>9$ generic points $p_i$ of \pr2 even when 
$n$ is not a square, is still open. 
(When $n$ is a square, it is easy to check that
$\epsilon(L, n)\ge1/\sqrt{n}-\varepsilon$ holds
for $n$ {\it general\/} points for 
any positive rational $\varepsilon$;
however, our generalization of this in \ir{introThm}
to any surface seems to be both new and nontrivial.)

Some of the best results obtained so far 
for the case that $l=1$ and $X=\pr2$ over ${\bf C}$ 
are due to Biran \cite{refBi}, who uses a powerful 
procedure for building
nef divisors. Although ad hoc applications of 
this procedure can yield impressive results
in particular cases (such as the calculation
$\epsilon(L, p_1,\ldots,p_{19})\ge 39/170$
in section 5 of \cite{refBi}, whereas our
\ir{introThm} gives only $39/171$), obtaining general
results by this procedure seems to require carefully 
constructed values of $n$. For example, given positive 
integers $a$ and $i$, Theorem 2.1A of \cite{refBi} 
gives bounds if $n=a^2i^2\pm2i$, or if 
$n=a^2i^2+i$ and $ai\ge3$. But in these cases 
there are certain positive integer solutions to $r^2-d^2n=1$:
for $n=a^2i^2\pm2i$, take $r=a^2i\pm1$ and $d=a$,
and for $n=a^2i^2+i$, take $r=2a^2i+1$ and $d=2a$;
either way Biran's bound is $\sqrt{1/n}\sqrt{1-{1/{r^2}}}$.
Applying \ir{introProp}(b), we recover as a 
special case of \ir{introThm} these same bounds
in those cases with $r\le n$, and we obtain even better bounds
via \ir{introPropB} when $n=a^2i^2\pm2i$ and $a=1$
(and hence $n+1$ is a square). 
For the cases when $i=1$ and either $n=a^2i^2-2i$ or 
$n=a^2i^2+i$, we have $r>n$ so \ir{introProp}(b)
does not apply, but (except in the case that
$n=a^2+1$ and either the characteristic is 2
or $a$ is a power of 2) \ir{introPropB} 
recovers Biran's bound via a refined application 
of our underlying approach; see \ir{rfnmnts}. 
Similarly, if $n=a^2i^2-i$, which \cite{refBi} does 
not treat (except in special cases when $n$ can also
be written in the form $n=a'^2i'^2\pm2i'$), 
we can take $d=2a$ and $r=2a^2i-1$
and again obtain the bound $\sqrt{1/n}\sqrt{1-{1/{r^2}}}$,
as long as $i>2$. (For bounds when $i\le2$, and more generally
when $n+2$ or $n+1$ is a square, see \ir{introPropB}.)

Another result over ${\bf C}$ for $X=\pr2$ and $l=1$ that should 
be mentioned is that $\epsilon(L, n)\ge 1/\sqrt{n+1}$ for $n\ge10$ 
general points \cite{refST}. Apart from cases
which follow from \cite{refnagprob} (which this paper generalizes)
and from those of \cite{refBi} just mentioned, 
this seems to have been the best estimate known up to now. However,
\ir{introThm} with \ir{introProp}(c) (or
\ir{introPropB} if $n\pm1$ is a square) is better in all cases.
Moreover, our approach applies to all surfaces in 
all characteristics. [Very recently, 
by a very elegant argument for $X=\pr2$ over ${\bf C}$, 
Szemberg \cite{refS} obtained a bound of the form
$\epsilon(L, n)\ge (1/\sqrt{n})\sqrt{1-1/(an)}$,
where $a$ currently can be as large as about 5.
But for $n$ sufficiently large, the 
bound of \ir{introThm} is, by \ir{introProp}(b)(iii), 
better except for a small fraction of cases. Nevertheless,
there are some special values of $n$
of particular interest, including certain small values
of $n$ and when $n=s^2-1$ where $s-1$ is a power of 2,
for which the bound of \cite{refS} is the best one we know.]

In short, Seshadri constants are difficult to compute
and in general remain unknown, but they are closely connected
to classical problems involving linear systems and thus are
of substantial interest. In this paper, using 
more broadly applicable geometric methods than
have been typical of work on this problem, we give a characteristic 
free approach to estimating Seshadri constants that nonetheless 
gives comprehensive improvements to currently known lower 
bounds. 

In preparation for stating our main result, let
$l$ and $n$ be positive integers and define the sets 
$$S_1(n,l)=\Bigl\{\;{r\over nd}\;\Bigl|\; 1\le r\le n, 
\quad1\le d,\quad {r\over d}\le \sqrt{nl}\;\Bigr\}$$
and 
$$S_2(n,l)=\Bigl\{\;{dl\over r}\;\Bigl|\; 1\le r\le n, \quad
1\le d,\quad{r\over d}\ge \sqrt{nl}\;\Bigr\}$$ of integer ratios.
Now define $S(n,l)=S_1(n,l)\cup S_2(n,l)$ and 
$\varepsilon_{n,l}=\hbox{max}(S(n,l))$. 
With a view to the important special case that 
$L$ is a line in $X=\pr2$, we will write
$\varepsilon_n$ for $\varepsilon_{n,1}$.

We now have the following result
(proved in \ir{scva} as \ir{scgen}):

\prclm{Theorem}{introThm}{Let $l=L^2$, where
$L$ is a very ample divisor on
an algebraic surface $X$. Then $\sqrt{l/n}\ge \epsilon(L, n)$,
and in addition, we have
$\epsilon(L, n)\ge\varepsilon_{n,l}$
unless $l\le n$ and $nl$ is a square,
in which case $\sqrt{l/n}=\varepsilon_{n,l}$
and $\epsilon(L, n)\ge\sqrt{l/n}-\varepsilon$
for every positive rational $\varepsilon$.}

\noindent(The somewhat awkward statement
in case $nl$ is a square is related to there possibly being
no {\it open} set of points such that $\epsilon(L, n)=\varepsilon_{n,l}$
in that case.)

Note that $\varepsilon_{n,l}$ is just 
the maximum element in the finite set
$$\Bigl\{\;{\lfloor d\sqrt{nl}\rfloor \over dn}\; \Bigl|\; 
1\le d\le \sqrt{{n\over l}}\;\Bigr\}
\cup\Bigl\{\;{1\over\lceil\sqrt{{n\over l}}\rceil}\;\Bigr\}\cup
\Bigl\{\;{dl\over\lceil d\sqrt{nl}\rceil}\; \Bigl| 
\;1\le d \le\sqrt{{n\over l}}\;\Bigr\}.$$
Thus for any given $n$ it is not hard to 
compute $\varepsilon_{n,l}$
exactly, even though it is not easy to give an explicit formula. 
As an alternative, we give some comparisons and
in addition determine $\varepsilon_{n,l}$ 
explicitly in some cases
(the proof of \ir{introProp} is in \ir{scva}):

\prclm{Proposition}{introProp}{Let $l$, $s$ and $n$ be positive integers.
\item{(a)} If $l\ge n$, then $\varepsilon_{n,l}=1$.
\item{(b)} Say $l<n$, let $d$ and $r\le n$ be positive integers
and put $\delta=r^2-nld^2$. 
\itemitem{(i)} If $nl$ is a square, then $\varepsilon_{n,l}=
\sqrt{l/n}$.
\itemitem{(ii)} We have
$$\varepsilon_{n,l}\ge\sqrt{l\over n}\sqrt{1-{\delta\over{r^2}}}
\hbox{ if $\delta\ge0$ and }
\varepsilon_{n,l}\ge\sqrt{l\over n}\sqrt{1+{\delta\over{nld^2}}}
\hbox{ if $\delta\le0$,}$$
with equality if $\delta=\pm1$.
\itemitem{(iii)} Moreover, 
$$\varepsilon_{n,l}>\sqrt{l\over n}\sqrt{1-{1\over n}}\eqno(*)$$
holds for at least half of the values of $l$ from $(n-1)/2$ to $n-1$
as long as $n>2$. Alternatively, given any positive integer $a$ and $s>2$, the
fraction of the number of values of $n$ in the range 
$s^2l\le n<(s+1)^2l$ for which
$$\varepsilon_{n,l}>\sqrt{l\over n}\sqrt{1-{1\over an}}\eqno(**)$$
fails to hold is at most 
$((2a^2-a+8)l+3)/(l(2s+1))$ if $a>2$, $(14l+3)/(l(2s+1))$ if $a=2$, 
and $(4l+2)/(l(2s+1))$ if $a=1$, and thus goes to 0 as $s$ increases.
\item{(c)} If $n\pm1$ is not a square, we have
$$\varepsilon_n>{1\over\sqrt{n+1}}>
\sqrt{{1\over n}}\sqrt{1-{1\over n}}.$$}

\vskip\baselineskip
For more explicit estimates of $\varepsilon_{n,l}$ 
and $\varepsilon_n$, see \ir{scspec} and \ir{ptwocor}. 
Also, for a given $n$, we note that the
first statement of \ir{introProp}(b)(iii)
significantly understates the number of $l$ from 1 to 
$n$ for which $(*)$ holds, which often is
$3n/4$ or more; see \ir{twthrds}.

Regarding \ir{introProp}(c), 
it is especially difficult to improve on previously known
bounds when $n$ is close to a square. If
$n\pm1$ or $n+2$ is a square, we can improve on $\varepsilon_n$ 
using a refinement of our basic approach. We obtain 
the following result, proved in \ir{rfnmnts}, which in 
all cases is better than $1/\sqrt{n+1}$, and recovers Biran's bound if
either $n+2$ is a square or (in certain cases) if $n-1$ is a square, 
and improves on Biran's bound if $n+1$ is a square:

\prclm{Proposition}{introPropB}{Let $n$ be a positive integer,
with $L$ a line in $X=\pr2$; we have:
\item{$\bullet$} if $9\le n+2$ is a square then 
$\epsilon(L,n)\ge \sqrt{1\over n}\sqrt{1-{1\over (n+1)^2}};$
\item{$\bullet$} if $9\le n+1$ is a square, then 
$\epsilon(L,n)\ge\sqrt{1\over n}\sqrt{1-{n-1\over n(\sqrt{n+1}+1)^2}},$
and (unless $\sqrt{n+1}-1$ is a power of 2 or the characteristic is 2)
$\epsilon(L,n)>\sqrt{1\over n}\sqrt{1-{1\over (\sqrt{n+1}-1)n}}$;
\item{$\bullet$} if $9\le n-1$ is a square, then 
$\epsilon(L,n)\ge \sqrt{1\over n}\sqrt{1-{n-1\over (n+\sqrt{n-1})^2}},$
and (unless $n-1$ is a power of 4 or the characteristic is 2)
$\epsilon(L,n)\ge \sqrt{1\over n}\sqrt{1-{1\over (2n-1)^2}}.$}

We include two corollaries that may be of interest.
For any $nl$ not a square, there are infinitely many
solutions $(r,d)$ to $r^2-d^2nl=1$. Unfortunately, 
if $r$ is too big we cannot apply \ir{introProp}(b)(ii)
to obtain a bound. By the next result (see
\ir{scva} for the proof), such solutions 
need not entirely go wasted:

\prclm{Corollary}{introCorA}{Let $L$ be a very 
ample divisor on a surface $X$ with $l=L^2$ and consider 
positive integers $n$, $r$ and $d=ab$.
If $r^2-nld^2=1$ and $a>b\sqrt{l/n}$, then 
$\pi^*L-(bl/r)(E_1+\cdots+E_{a^2n})$ is a nef 
${\bf Q}$-divisor,
where $\pi:Y\to X$ is the birational morphism obtained
by blowing up $a^2n$ general points $p_i$, $E_i$
being the exceptional divisor corresponding to $p_i$;
in particular, we have $\epsilon(L, a^2n)\ge bl/r$.}

The preceding result is suggestive of the procedure
of \cite{refBi} on $X=\pr2$ over ${\bf C}$, which, 
for example, can easily be used to show that 
$at\pi^*L-m(E_1+\cdots+E_{a^2n})$ is nef if 
$t\pi^*L-m(E_1+\cdots+E_{n})$ is. Similarly, the next 
corollary (proved in \ir{scva}) generalizes two additional 
facts known on $X=\pr2$ over ${\bf C}$: that the divisor
$H=t\pi^*L-(E_1+\cdots+E_{n})$ is ample if
$H^2>0$ (see \cite{refXa} or \cite{refKub})
and $H=t\pi^*L-2(E_1+\cdots+E_{n})$ is nef
if $H^2\ge0$ (see Theorem 2.1.B of \cite{refBi}).

\prclm{Corollary}{introCorB}{Let $L$ be a very 
ample divisor on a surface $X$ with $l=L^2$ and consider 
positive integers $n>l$ and $d\le \sqrt{n/l}$.
Let $\pi:Y\to X$ be the birational morphism obtained
by blowing up $n$ general points $p_i$, $E_i$
being the corresponding exceptional divisor.
Then $H_r=r(\pi^*L)-dl(E_1+\cdots+E_n)$ is 
a nef divisor for all integers $r>d\sqrt{nl}$,
and an ample ${\bf Q}$-divisor for all rationals
$r>\lceil d\sqrt{nl}\rceil$.}
 
Our approach uses an explicit construction 
in \ir{nefblups} of nef divisors  
on the blow up $Y$ of $X$ at the points $p_i$, with the
points taken in special position. As a consequence of these nef divisors,
we obtain various bounds in \ir{scva}. 
The construction of nef divisors given in
\ir{nefblups} can be refined to sometimes obtain better bounds.
Since doing this introduces some complications, we
segregate this material to \ir{rfnmnts}. In \ir{appls} we
discuss additional applications of the existence of 
these nef divisors to various classical problems involving 
linear systems on \pr2. Analogous remarks could be made
for surfaces more generally, but complications 
(such as irregular surfaces, and failure of 
vanishing theorems to hold in certain circumstances
in positive characteristics) arise that would require 
special treatment. Thus we leave such remarks for 
the reader to work out in cases of his or her 
own interest.

\irrnSection{Nef Divisors on Blow ups}{nefblups}
The foundation for our
results is a method for constructing nef divisors generalizing
what is done in \cite{refnagprob}. To state our basic lemma,
let $X_0$ be an algebraic surface. We will call a sequence 
$p_1,\cdots,p_n$ of points a {\it proximity sequence\/} if
$p_1\in X_0$ is smooth, $p_2$ is a point of the exceptional divisor
of the blow up $\pi_1: X_1\to X_0$ at $p_1$, and, for 
$2<i<n$, $\pi_{i-1}:X_{i-1}\to X_{i-2}$ is the blow up of 
$X_{i-2}$ at $p_{i-1}$ and $p_i$ is a point of the exceptional 
divisor of $\pi_{i-1}$, but not a point of the proper transform
of the exceptional divisor of $\pi_{i-2}$. Denote the composition
$\pi_n\circ\cdots\circ\pi_1$ by $\pi: X_n\to X_0$ and denote 
by $E_i$ the scheme theoretic fiber 
$(\pi_n\circ\cdots\circ\pi_{i-1})^{-1}(p_i)$. Thus $E_i$ is 
a divisor on $X_n$, the total transform of $p_i$, and, for each 
$1\le i<n$, $[E_i-E_{i+1}]$ is the class of a reduced irreducible
divisor. 

\prclm{Lemma}{neflemA}{Let $\pi:Y\to X$ be the morphism
obtained by blowing up general points $p_1,\cdots,p_n$
of an algebraic surface $X$, and let $E_i$ be the exceptional 
divisor corresponding to each point $p_i$.
Let $\pi':Y'\to X$ be the morphism corresponding to
a proximity sequence $p'_1,\cdots,p'_n$ for $X$, with $E'_i$ 
being the exceptional divisor on $Y'$ corresponding to $p'_i$.
Suppose we are given a divisor $L$ on $X$ with $l=L^2$,
and integers $1\le d$, $1\le r\le n$, and 
$m_1\ge\cdots\ge m_r\ge 0$ such that 
$[d(\pi'^*L)-m_1E'_1-\cdots-m_rE'_r]$
is the class of an irreducible effective divisor $C$.
Then $a_0d(\pi^*L)-a_1E_1-\cdots-a_nE_n$
is a nef ${\bf Q}$-divisor on $Y$ for any rational numbers $a_i$
satisfying:
\item{$\bullet$} $a_1\ge a_2\ge \cdots\ge a_n\ge 0$;
\item{$\bullet$} $a_0d^2l\ge a_1m_1+\cdots+a_rm_r$;
\item{$\bullet$} $a_0^2d^2l>a_1^2+\cdots+a_n^2$;
\item{$\bullet$} $(m_1+\cdots+m_i)a_0\ge a_1+\cdots+a_i$
for $1\le i\le r$; and
\item{$\bullet$} $(m_1+\cdots+m_r)a_0\ge a_1+\cdots+a_n$.}

\noindent{\bf Proof}: The point is to show
$N=a_0d(\pi'^*L)-a_1E'_1-\cdots-a_nE'_n$ is nef on $Y'$.
Since by assumption $N^2>0$, 
the result then follows since being nef and of 
positive self-intersection is an open condition
for flat families of line bundles. But the divisor
class $[N]$ is nef on $Y'$ because it is a nonnegative ${\bf Q}$-linear 
combination of classes of irreducible effective divisors, 
each of which $N$ meets nonnegatively. In particular, the 
last two bulleted hypotheses guarantee that $N$ is the sum of 
$a_0C$ with various nonnegative multiples of
$E_i-E_{i+1}$ for $1\le i<n$ and $E_n$. The first two
bulleted hypotheses guarantee that $N$ meets each summand
nonnegatively. (It follows that $a_0^2d^2l\ge a_1^2+\cdots+a_n^2$,
i.e., $N^2\ge0$, is a consequence of the other hypotheses,
and that $N^2>0$ is automatic unless:
$a_0d^2l= a_1m_1+\cdots+a_rm_r$, and 
$((m_1+\cdots+m_i)a_0-(a_1+\cdots+a_i))(a_i-a_{i+1})=0$
for $1\le i<r$, and $(a_r-a_i)a_i=0$ for $i>r$, and 
$(m_1+\cdots+m_r)a_0=a_1+\cdots+a_n$.)
\qed

We now apply the preceding lemma in case $L$ is a very ample 
divisor on $X$.

\prclm{Lemma}{neflem}{Let $\pi:Y\to X$ be the morphism 
blowing up general points $p_1,\ldots,p_n$ of
an algebraic surface $X$, let $E_1,\cdots,E_n$ be the 
exceptional divisors corresponding to the points $p_i$,
let $L$ be a very ample divisor on $X$ and put $l=L^2$
and $L'=\pi^*L$. 
Given positive integers $r\le n$ and $d$,
and nonnegative rational numbers (not all 0)
$a_0\ge a_1\ge \cdots \ge a_n\ge 0$
such that $a_0d^2l\ge a_1+\cdots+a_r$,
$a_0^2d^2l> a_1^2+\cdots+a_n^2$
and $ra_0\ge a_1+\cdots+a_n$, then
$(a_0d)L'-(a_1E_1+\cdots+a_nE_n)$ is a 
nef ${\bf Q}$-divisor on $Y$.}

\noindent{\bf Proof}: Since $L$ is very ample,
$|dL|$ has an irreducible member $C'$ passing through 
and smooth at some smooth point $p'_1\in X$.  
Take as our proximity sequence points of $C'$ infinitely 
near to $p'_1$; i.e., blow up $p'_1$ and take
$p'_2$ to be the point of the proper transform of $C'$
infinitely near to $p'_1$. Similarly define $p'_i$ for
all $i\le r$. Then extend to a proximity sequence
$p'_1,\ldots,p'_n$ such that $p'_{r+1}$ is {\it not\/}
on the proper transform of $C'$. Blowing up the points
of the sequence gives the morphism $\pi':Y'\to X$,
and $[d(\pi'^*L)-E'_1-\cdots-E'_r]$ is the class
of the proper transform of $C'$, which is irreducible.
It is now easy to check that the hypotheses of
\ir{neflemA} apply (take $m_i=1$ for all $i$),
giving the result.\qed

One application of \ir{neflem} is
to provide nef divisors $F$ which can be employed to test 
for effectivity: given a divisor
$H$ on $X$ and integers $b_i$ (we may as well assume 
$b_1\ge \cdots \ge b_n\ge 0$), if $F\cdot (\pi^*H-b_1E_1-\cdots-b_nE_n)<0$
for some nef $F$, then $|\pi^*H-b_1E_1-\cdots-b_nE_n|$ is empty. 
Given $\pi^*H-b_1E_1-\cdots-b_nE_n$, the optimal nef test divisor
$F$ provided by \ir{neflem} can be found by
linear programming. (Keep in mind that we can always 
normalize so that $a_0=1$, and clearly one need consider only
finitely many $r$ and $d$.)

In order to avoid linear programming,
the following corollary obtains
some special cases of particular interest.

\prclm{Corollary}{nefcor}{Given $X$, $Y$, $l$, $n$ 
and $L'$ as in \ir{neflem}, let $d$ and $r\le n$ 
be positive integers. Then we have the following cases.
\item{(a)} If $r^2>nd^2l$, then 
$rL'-dl(E_1+\cdots+E_n)$ is nef.
\item{(b)} If $r^2<nd^2l$, then 
$ndL'-r(E_1+\cdots+E_n)$ is nef.
\item{(c)} If $r^2=nd^2l$, then 
$tL'-(E_1+\cdots+E_n)$ is nef for all rationals
$t>\sqrt{n/l}$.}

\noindent{\bf Proof}: Apply \ir{neflem} for various values of the $a_i$.
For (a), take $a_0=r/d$ and $a_1=\cdots=a_n=ld$.
For (b), take $a_0=n$ and $a_i=r$, $i>0$.
For (c), take $a_0>n/r$ and $a_i=1$, $i>0$.
\qed

There are cases when one may want to construct nonuniform
nef divisors. Here are some examples of such.

\prclm{Corollary}{nefcorB}{Given $X$, $Y$, $l$, $n$ 
and $L'$ as in \ir{neflem},
let $d$ and $r\le n$ be positive integers.
Then we have the following cases.
\item{(a)} If $d^2l>r$, then $dL'-(E_1+\cdots+E_r)$ is nef.
\item{(b)} If $d^2l\le r$, then $d'L'-(E_1+\cdots+E_{ld^2})$ is a
nef ${\bf Q}$-divisor for all rational $d'>d$,
and, for each integer $1\le j< d^2l$, 
$$d'L'-(E_1+\cdots+E_j)-{(d^2l-j)\over(r-j)}(E_{j+1}+
\cdots+E_{\lfloor \lambda\rfloor}+
(\lambda-\lfloor\lambda\rfloor)E_{\lceil\lambda\rceil})$$
is a nef ${\bf Q}$-divisor, where 
$\lambda=\hbox{min}\{r+(r-d^2l)(r-j)/(d^2l-j),n\}$
and $d'\ge d$ is any rational such that 
$d'>d$ if $\lambda=\lceil\lambda\rceil\le n$.}

\noindent{\bf Proof}: We apply \ir{neflem}.
For (a), take $a_i=1$ for $0\le i\le r$ and $a_i=0$ for $i>r$.
For the first part of (b), take $a_0=d'/d$ and $a_i=1$ for 
$0< i\le d^2l$ and $a_i=0$ for $i>d^2l$.
For the rest, the idea is to choose $a_i$ such that
$a_i=1$ for $0\le i\le j$, with the $a_i$ for $j<i\le r$
being equal and as large as possible subject to
$a_0d^2l\ge a_1+\cdots+a_r$ (hence $a_i=(d^2l-j)/(r-j)$
for $j<i\le r$), and finally for as many as possible of the
remaining $a_i$ also to equal $(d^2l-j)/(r-j)$, subject to
$ra_0\ge a_1+\cdots+a_n$. Thus we take $a_i=(d^2l-j)/(r-j)$
for $r<i\le\lfloor\lambda\rfloor$, $a_i=0$ for $i>\lceil\lambda\rceil$
and, if $(\lambda-\lfloor\lambda\rfloor)>0$, we take 
$a_{\lceil\lambda\rceil}=(\lambda-\lfloor\lambda\rfloor)(d^2l-j)/(r-j)$
(in which case $ra_0\ge a_1+\cdots+a_n$ will be an equality).
The requirement on $d'$ ensures positive self-intersection.
\qed

\irrnSection{Seshadri Constants of Very Ample Divisors}{scva}
For a very ample divisor $L$ on a surface $X$ with $L^2=l$
and any $n$ distinct points $p_i$ on $X$, it is easy to see that 
$\epsilon(L,p_1,\ldots,p_n)\le\sqrt{l/n}$: just note that
for any $\varepsilon$ bigger than $\sqrt{l/n}$ we can find a rational
$\delta<\sqrt{l/n}$ such that $F_\varepsilon\cdot F_\delta<0$,
where $F_t=L'-t(E_1+\cdots+E_n)$. But $F_\delta^2>0$, so
for appropriate integers $N$ sufficiently large, $|NF_\delta|$
is nonempty, hence $\epsilon(L,p_1,\ldots,p_n)\le\varepsilon$.

Lower bounds are more difficult.
By applying \ir{nefcor}, we establish
our main lower bound.

\prclm{Theorem}{scgen}{Let $l=L^2$, where
$L$ is a very ample divisor on
an algebraic surface $X$. Then $\sqrt{l/n}\ge \epsilon(L, n)$,
and in addition, we have
$\epsilon(L, n)\ge\varepsilon_{n,l}$
unless $l\le n$ and $nl$ is a square,
in which case $\sqrt{l/n}=\varepsilon_{n,l}$
and $\epsilon(L, n)\ge\sqrt{l/n}-\varepsilon$
for every positive rational $\varepsilon$.}

\noindent{\bf Proof}: We noted $\sqrt{l/n}\ge \epsilon(L, n)$ above. 
If $l>n$, then $\varepsilon_{n,l}=1$ by \ir{introProp},
but $L'-E_1-\cdots-E_n$ is nef by \ir{nefcor}(b)
(take $r=n$ and $d=1$), so $\epsilon(L, n)\ge1$.
If $l<n$ but $nl$ is not a square, then 
$\epsilon(L, n)\ge\varepsilon_{n,l}$ follows from 
\ir{nefcor}, parts (a) and (b). Finally, if $l\le n$ and $nl$ 
is a square, then $\sqrt{l/n}=\varepsilon_{n,l}$
follows from \ir{introProp}, and 
$\epsilon(L, n)\ge\sqrt{l/n}-\varepsilon$
holds for every positive rational $\varepsilon$
by \ir{nefcor}(c).\qed

Although one needs to check only finitely many values of 
$r$ and $d$ to compute $\varepsilon_{n,l}$, it is nonetheless
useful to have more explicit lower bounds. For that purpose, 
given positive integers $n$ and $l$, let $d^*=\lceil\sqrt{n/l}\rceil$,
$d_*=\lfloor\sqrt{n/l}\rfloor$,
$r^*=\lceil d_*\sqrt{nl}\rceil$, and
$r_*=\lfloor d_*\sqrt{nl}\rfloor$.

\prclm{Corollary}{scspec}{Let $l$ and $n$ be positive integers.
Then $\varepsilon_{n,l}\ge 1/d^*$, and, if
$l\le n$, then also $\varepsilon_{n,l}\ge 
\hbox{max}(r_*/(nd_*), d_*l/r^*)$.}

\noindent{\bf Proof}: For the first inequality, use $r=n$ and $d=d^*$,
and check that then $1\le r\le n$, $1\le d$, and $r/d\le \sqrt{nl}$,
so in this case $r/(nd)\in S_1$.
For $r_*/(nd_*)$ in the second inequality, use $r=r_*$ and $d=d_*$,
and again check that $1\le r\le n$, $1\le d$ (because $l\le n$), 
and $r/d\le \sqrt{nl}$, so $r/(nd)\in S_1$.
For $d_*l/r^*$, use $r=r^*$ and $d=d_*$,
and check that $1\le r\le n$, $1\le d$, 
and $r/d\ge \sqrt{nl}$, so $dl/r\in S_2$. \qed

\rem{Remark}{optrd} The values of $r$ and $d$ obtained
using $d^*$, $d_*$, $r^*$, and $r_*$ are not always optimal.
For example, if $n=33$ and $l=1$, then 
$\varepsilon_{n,l}=4/23\in S_2(n,l)$ comes from $r=23$ and $d=4$,
but $1/d^*=1/6$, $r_*/(nd_*)=28/(33\cdot5)$, $d_*/r^*=5/29$
are all less than $4/23$. 

\rem{Remark}{ptworem} If $L$ is a line in $X=\pr2$,
then $l=1$. If we denote $\lfloor\sqrt{n}\rfloor$ by $s$
and $\lfloor(n-s^2)/2\rfloor$ by $t$, then either 
$n=s^2+2t$ or $n=s^2+2t+1$, where $0\le t\le s$
(with $t<s$ in the latter case). With respect to $s$ and $t$
in the case that $n$ is not a perfect square,
it is not hard to check that $r^*=s^2+t$, $r_*=s^2+t-1$, $d^*=s+1$ 
and $d_*=s$ if $n=s^2+2t$, while $r^*=s^2+t+1$, $r_*=s^2+t$, $d^*=s+1$ 
and $d_*=s$ if $n=s^2+2t+1$. If $n=s^2$, then
$r^*=r_*=s^2$ and $d^*=d_*=s$. (For a more symmetrical treatment,
under some restrictions, of cases (a) and (b) of the 
following corollary, see \ir{nminusoneexmpl}.)

\prclm{Corollary}{ptwocor}{Let $1\le s$ and $0\le t\le s$ be integers.
\item{(a)} If $n=s^2+2t$, then $\varepsilon_{n}\ge s/(s^2+t)$.
\item{(b)} If $n=s^2+2t+1$ and $t<s$, then 
$\varepsilon_{n}\ge (s^2+t)/(s(s^2+2t+1))$
\item{(c)} If $s>1$ and $n=s^2+2t+1$ and $0<t<(\sqrt{2}-1)(s-1)$, then 
$\varepsilon_{n}\ge (s(s-1)+t)/((s-1)(s^2+2t+1))>(s^2+t)/(s(s^2+2t+1))$.
\item{(d)} If $s>1$ and $n=s^2+2t+1$ and 
$(\sqrt{2}-1)(s-1)<t<\sqrt{1.25s^2-s}-s/2$, then 
$\varepsilon_{n}\ge (s-1)/(s(s-1)+t)>(s^2+t)/(s(s^2+2t+1))$.}

\noindent{\bf Proof}: For (a) and (b), apply \ir{scspec}, 
using the expressions for $r^*$, $r_*$, $d^*$ and $d_*$ 
in \ir{ptworem}. For (c) and (d), apply \ir{nefcor}(b) and (a), resp., 
with $r=s(s-1)+t$ and $d=s-1$.\qed

\rem{Remark}{Xurem} For $L$ a line in $X=\pr2$ over ${\bf C}$, 
\cite{refST} proves for $n\ge10$ that $\epsilon(L,n)\ge1/\sqrt{n+1}$. 
It is easy to check that \ir{ptwocor} gives 
a better result in all cases except when $n\pm1$ is a square.
For improvements in these cases, see \ir{introPropB}.
\vskip\baselineskip

\rem{Remark}{smalln} The lower bounds given in \ir{ptwocor}(a,b)
actually equal both $\epsilon(L, n)$ and
$\varepsilon_{n}$ if $1\le n\le6$. For $n=7$, 
$\epsilon(L, 7)=3/8$ (since
$F=8L'-3E_1-\cdots-3E_7$ is known to be nef while
$E=3L'-2E_1-E_2-\cdots-E_7$ is effective with $F\cdot E=0$),
and $\epsilon(L, 8)=6/17$ for $n=8$ (since
$F=17L'-6E_1-\cdots-6E_8$ is nef while
$E=6L'-3E_1-2E_2-\cdots-2E_8$ is effective with $F\cdot E=0$),
whereas in fact
$\varepsilon_{7}=5/14$ and $\varepsilon_{8}=1/3$.
\vskip\baselineskip

\noindent{\bf Proof} of \ir{introProp}: 
Part (c) is easy to check using \ir{ptwocor}(a,b).
For part (a), note that we have
$\varepsilon_{n,l}\ge 1/d^*$ by \ir{scspec}, but $d^*=1$
for $l\ge n$. On the other hand, by definition either
$\varepsilon_{n,l}=r/(nd)$ for some positive $r$ and $d$
with $r\le n$ (in which case clearly $\varepsilon_{n,l}\le1$),
or $\varepsilon_{n,l}=dl/r$ for some positive $r$ and $d$
with $r\le n$ and $r^2\ge d^2nl$ 
(hence $dl/r=d^2nl/(rdn)\le r^2/(rdn)\le 1$).

Consider part (b). Given $r$ and $d$ with 
$\delta=r^2-d^2nl$, it is easy to check that
$dl/r=\sqrt{l/n}\sqrt{1-\delta/r^2}$ if $0\le\delta$, while 
$r/(nd)=\sqrt{l/n}\sqrt{1+\delta/(d^2nl)}$ if $\delta\le0$.
The inequalities in (b)(ii) now follow by definition of $\varepsilon_{n,l}$.
Moreover, this makes it clear that $\varepsilon_{n,l}\le\sqrt{l/n}$,
so if $nl=q^2$ for some $q$, we can take
$r=q$ and $d=1$ to see $\varepsilon_{n,l}\ge r/(nd)=\sqrt{l/n}$.
This proves part (b)(i). To prove the statement about
equality in (b)(ii), first assume $r^2-nld^2=1$
with $r\le n$. It suffices to
show $\varepsilon_{n,l}=dl/r$.

For any positive integer $t\le \sqrt{n/l}$, 
denote $\lceil t\sqrt{nl}\rceil$ by $r_t$;
e.g., we have $r=r_d$.
Since $r^2-nld^2=1$, we know that $nl$ is not a perfect square,
so $\lfloor t\sqrt{nl}\rfloor=r_t-1$. 
Now, $\varepsilon_{n,l}$ is just the maximum in 
$\{(r_t-1)/(tn) | 1\le t\le \sqrt{n/l}\}\cup
\{1/t | t=\lceil\sqrt{n/l}\rceil\}
\cup \{tl/r_t | 1\le t\le \sqrt{n/l}\}$.
We will show that $r_t=\lceil rt/d\rceil$. Assuming this,
it follows that $dl/r=tl/(rt/d)\ge tl/r_t$, so $dl/r$
is the maximum of $\{tl/r_t | 1\le t\le \sqrt{n/l}\}$.
We must also show $dl/r$ is as large as every element of
$\{(r_t-1)/(tn) | 1\le t\le \sqrt{n/l}\}\cup
\{1/\lceil\sqrt{n/l}\rceil\}$.
But from $r^2-nld^2=1$ we derive $r^2/(d^2l^2)=n/l+1/(d^2l^2)$.
If $dl/r<1/\lceil\sqrt{n/l}\rceil$, then 
$\lceil\sqrt{n/l}\rceil<r/(dl)$, and there must be an integer
$k$ with  $\sqrt{n/l}\le k<r/(dl)$, hence
$r^2/(d^2l^2)-1/(d^2l^2)=n/l\le k^2<r^2/(dl)^2$,
and so $r^2-1\le k^2d^2l^2<r^2$, which is absurd.
As for $\{(r_t-1)/(tn) | 1\le t\le \sqrt{n/l}\}$,
we have $rt=r_td-\rho$ where $0\le \rho<d$. By solving for $r_t$
and substituting, we see $(r_t-1)/(nt)\le dl/r$ 
if and only if $t=t(r^2-nld^2) \le (d-\rho)r$. But $(d-\rho)r
\ge r>d\sqrt{nl}\ge\sqrt{nl}\ge \sqrt{n/l}\ge t$, as required.

We are left with checking $r_t=\lceil rt/d\rceil$. 
Since $r=\lceil d\sqrt{nl}\rceil$,
we see $rt\ge dt\sqrt{nl}$, so $r_t=\lceil rt/d\rceil$ follows 
if we show there is no integer $k$ with $t\sqrt{nl}<k<rt/d$ (equivalently,
that there is no $k$ with $t^2d^2nl<k^2d^2<r^2t^2$),
but such a $k$ would imply $t^2d^2nl<(rt-1)^2$. Thus it suffices
to show that  $rt - dt\sqrt{nl}\le1$ for $1\le t\le \sqrt{n/l}$.
But $t(r-d\sqrt{nl})=t(r^2-nld^2)/(r+d\sqrt{nl})<t/(2d\sqrt{nl}) 
\le \sqrt{n/l}/(2d\sqrt{nl})=1/(2dl)\le 1$.

Suppose now that $r^2-nld^2=-1$, in which case
we must show $\varepsilon_{n,l}=r/(nd)$.
This time $r_t=\lfloor t\sqrt{nl}\rfloor$,
$\lceil t\sqrt{nl}\rceil=r_t+1$ and we will show
$r_t=\lfloor rt/d\rfloor$. 
It follows that $r/(nd)=rt/(ndt)\ge r_t/(nt)$.
We also have $r/(nd)\ge 1/\lceil\sqrt{n/l}\rceil$:
if not then $\lceil\sqrt{n/l}\rceil<nd/r$,
but $n/l=(dn/r)^2-n/(r^2l)$, so
there is an integer $k$ with
$(dn/r)^2-n/(r^2l)=n/l\le k^2<(dn/r)^2$,
hence $(dn)^2-n/l\le k^2r^2<(dn)^2$,
but this is not sufficient distance between squares unless
$d=n=l=1$, which contradicts $r^2-d^2nl=-1$.
Now compare $r/(nd)$ with $tl/(r_t+1)$.
We have $rt=r_td+\rho$ where $0\le \rho<d$,
and as before $tl/(r_t+1)\le r/(nd)$ 
if and only if $t\le (d-\rho)r$. But $(d-\rho)r
\ge r>d\sqrt{nl}-1=dl\sqrt{n/l}-1$,
and $dl\sqrt{n/l}-1\ge t$
unless $t=\lfloor\sqrt{n/l}\rfloor$ and $dl=1$,
but in that case it is easy to check that $t=r-1$.

We are left with checking $r_t=\lfloor rt/d\rfloor$. 
Since $r=\lfloor d\sqrt{nl}\rfloor$,
we see $rt\le dt\sqrt{nl}$, so $r_t=\lfloor rt/d\rfloor$ follows 
if $dt\sqrt{nl}-rt\le1$ for $1\le t\le \sqrt{n/l}$.
But $t(d\sqrt{nl}-r)=t(nld^2-r^2)/(d\sqrt{nl}+r)<t/(d\sqrt{nl}) 
\le \sqrt{n/l}/(dl\sqrt{n/l})\le1/(dl)\le 1$.

Finally, consider part (b)(iii). There exist $r$ and $d$ such that either
$\varepsilon_{n,l}=dl/r$ with $0\le\delta=r^2-d^2nl$ or
$\varepsilon_{n,l}= r/(nd)$ with $\delta\le0$.
If $0\le\delta$, it's enough to check
that $dl/r>\sqrt{l/n}\sqrt{1-1/n}$, but as above
$dl/r=\sqrt{l/n}\sqrt{1-\delta/r^2}$ so it suffices to check that 
$\delta/r^2<1/n$; i.e., that $\delta<r^2/n$.
If $\delta\le0$, the argument is the same except
$\varepsilon_{n,l}=r/(nd)=\sqrt{l/n}\sqrt{1+\delta/(d^2nl)}$.

So, to bound the number of $l$ for which $(*)$ in (b)(iii)
holds, we check whether $-d^2l<\delta< r^2/n$ holds
when $d=1$ and $r$ is either $r=\lfloor\sqrt{nl}\rfloor$ or
$r=\lceil\sqrt{nl}\rceil$. But $-d^2l<\delta< r^2/n$ holds
if either $\lceil\sqrt{nl}\rceil< \sqrt{nl}/\sqrt{1-1/n}$
or $\lfloor\sqrt{nl}\rfloor> \sqrt{nl}\sqrt{1-1/n}$, which
is equivalent to having the interval 
$I_l=(\sqrt{nl}\sqrt{1-1/n},\sqrt{nl}/\sqrt{1-1/n})$
contain an integer. 
It is not too hard to check that the union 
$I_{\lceil(n-1)/2\rceil}\cup\cdots\cup I_{n-1}$
contains the interval $(\sqrt{n(n-1)/2},n-1)$,
and thus the number of values of $l$ between $(n-1)/2$ and $n-1$
for which $\varepsilon_{n,l}>\sqrt{l/n}(\sqrt{1-1/n})$ holds
is always at least $(n-1)-\sqrt{n(n-1)/2}-1$. This is at least half
of the number of $l$ in the range $(n-1)/2\le l <n$, as long as
$n\ge45$. An explicit check for $3\le n\le 44$ shows that
$(*)$ still holds for at least half
of the number of $l$ in the range $(n-1)/2\le l <n$.

For $(**)$, we apply Dirichlet's theorem from elementary number theory,
which says there are integers $0<r<n+1$ and $d\ge1$ such that
$|r/\sqrt{nl}-d|\le 1/(n+1)$. Given such an $r$ and $d$,
we have $|r-d\sqrt{nl}|\le \sqrt{nl}/(n+1)$, hence
$|\delta|=|r^2-nld^2|\le (r+d\sqrt{nl})\sqrt{nl}/(n+1)$.
Thus, if $\delta<0$, we have $|\delta|/(nld^2) < 
\sqrt{nl}(r+d\sqrt{nl})/(n^2ld^2)<\sqrt{nl}(2d\sqrt{nl})/(n^2ld^2)
=2/(nd)$, and if $\delta>0$, we have $\delta/r^2 < 
\sqrt{nl}(r+d\sqrt{nl})/(nr^2)<\sqrt{nl}(2r)/(nr^2)<(r/d)(2r)/(nr^2)
=2/(nd)$. 

For some $r=r'$ and $d=d'$ and $\delta=r'^2-nld'^2$, we know that 
$\varepsilon_{n,l}=\sqrt{l/n}\sqrt{1-x}$, where $x=|\delta|/(nld'^2)$
if $\delta<0$ and $x=\delta/r'^2$ if $\delta>0$, and thus 
either way $x<2/(nd')$. It follows that if, for the 
given $n$, $|r/\sqrt{nl}-d|\le 1/(n+1)$
holds for no $r$ and $d$ with $0<r<n+1$ and $1\le d<2a$, then 
$x<2/(nd')\le 1/(an)$, and hence
$(**)$ holds for this $n$. So to count
those $n$ in the range $s^2l\le n<(s+1)^2l$ for which
$(**)$ holds, it is enough to count how often
$|r/\sqrt{nl}-d|\le 1/(n+1)$ holds for $1\le d<2a$ and $0<r<n+1$.

But for any given $d$, 
$|r/\sqrt{nl}-d|\le 1/(n+1)$ holds for some $r$ only if
the interval $(d\sqrt{nl}-\sqrt{l/n},d\sqrt{nl}+\sqrt{l/n})$
contains an integer, i.e., only if
$[d^2nl-2dl+1,d^2nl+2dl]$ contains a square.
Now, the interval $[d^2s^2l^2,d^2(s+1)^2l^2]$ 
contains $dl+1$ squares, and we are interested 
in counting for how many $n\in [s^2l,(s+1)^2l)$
does the interval $[d^2nl-2dl+1,d^2nl+2dl]$ 
contain one of these squares. Since $s>2$, no 
interval $[d^2nl-2dl+1,d^2nl+2dl]$ can contain two squares,
and for $d>3$, the intervals are disjoint and $d^2(s+1)^2l^2$
is in no interval, so at most $dl$ of the intervals 
contain squares. For $2\le d\le3$,
consecutive intervals overlap but no point 
lies in three intervals (and $d^2(s+1)^2l^2$
is in no interval if $d=3$ and only in 
the last interval if $d=2$),
so there are at most $2dl$ intervals that contain squares
when $d=3$ and at most $2dl+1$ when $d=2$.
Similarly, for $d=1$ at most four intervals 
overlap at a single point and $d^2(s+1)^2l^2$
is in two intervals,
so there are at most $4dl+2$ intervals that contain squares.
(These are of course typically overestimates since 
some squares may lie in no intervals.) Summing over $1\le d<2a$,
we find that of the $n$ in the range $s^2l\le n<(s+1)^2l$
there are at most $(4l+2+4l+1+6l)+(4l+5l+\cdots+(2a-1)l)$
(i.e., $(2a^2-a+8)l+3$ if $a>2$, $14l+3$ if $a=2$, 
and $4l+2$ if $a=1$) values of $n$ whose corresponding
interval contains a square. \qed

\rem{Remark}{twthrds} Our estimate that $(*)$ in \ir{introProp}(b)
holds for at least half of $(n-1)/2\le l<n$ 
understates how often $(*)$ holds. One reason for this is that
the intervals $I_l$ in the proof overlap, and thus
the same integer can lie in more than one interval,
but our estimate counts only some of those integers,
and at most once each. Also, our estimate is based on
a check only for $d=1$. We can partially account
for cases with $d>1$ by again applying
Dirichlet's theorem, as in the proof of \ir{introProp}(b)(iii).
As we saw there, $-2dl<\delta< 2r^2/(nd)$ holds for any $r$ and $d$
such that $|r/\sqrt{nl}-d|\le 1/(n+1)$ with $0<r<n+1$ and $d\ge1$. 
Therefore, $-d^2l<\delta< r^2/n$ also holds
if in addition $d>1$. Thus, as long as $|r/\sqrt{nl}-1|>1/(n+1)$
holds for $r=\lfloor\sqrt{nl}\rfloor$ and $r=\lceil\sqrt{nl}\rceil$
(which we can rewrite as $(n+2)\sqrt{nl}/(n+1) -1<
\lfloor\sqrt{nl}\rfloor<n\sqrt{nl}/(n+1)$),
we see that the solution to $|r/\sqrt{nl}-d|\le 1/(n+1)$
guaranteed by Dirichlet's theorem 
must have $d>1$ and hence $(*)$ holds. 
I.e., if the interval $J_l=((n+2)\sqrt{nl}/(n+1) -1,n\sqrt{nl}/(n+1))$
contains an integer, then $\varepsilon_{n,l}>\sqrt{l/n}\sqrt{1-1/n}$.

One can check that the intervals $J_l$ are nonempty
as long as $l$ is less than about $n/4$, and that
these intervals and the $I_l$ are all disjoint  
as long as $l$ is less than about $n/2$, and that
the union of the $I_l$ for $l$ more than about $n/2$
is about $(n/\sqrt{2},n)$. Thus the union of all of the intervals
$I_l$ and $J_l$ has measure about $0.61n$, 
so it is reasonable (but not guaranteed) to expect that
at least $61\%$ of the values of $l$ from 1 to $n$ should
give $\varepsilon_{n,l}>\sqrt{l/n}\sqrt{1-1/n}$. 
To take into account overlaps among the $I_l$ for $l>n/2$,
we might instead want to consider the sum of the lengths of 
the intervals. This is about $3n/4$, and so it is reasonable to expect
that typically at least $75\%$ of the values of $l$
from 1 to $n$ result in $\varepsilon_{n,l}>\sqrt{l/n}\sqrt{1-1/n}$.
Explicit computations for various $n$ show, in fact, that
percentages around $80\%$ are common. For $n$ from 15 to 200, 
the smallest percentage ($63\%$) occurs for $n=19$ and the largest
($87.6\%$) for $n=97$. For some larger $n$, we have 
$85\%$ for $n=313$, $75\%$ for $n=314$, $78\%$ for $n=3079$
and $80.8\%$ for $n=3080$.
\vskip\baselineskip

We now prove \ir{introCorA} and \ir{introCorB}.
\vskip\baselineskip

\noindent{\bf Proof} of \ir{introCorA}: This follows from
\ir{nefcor}(a), if we check that $r^2>b^2(a^2n)l$
and $r\le a^2n$. But $r^2=1+d^2nl>d^2nl=(ab)^2nl$,
and $r^2-1=d^2nl=a^2nb^2l=a^2n^2b^2l/n<a^4n^2$,
hence $r^2\le a^4n^2$, as required.
\qed

\noindent{\bf Proof} of \ir{introCorB}: For the first part it is
clearly enough to consider the case that $r=\lceil d\sqrt{nl}\rceil$,
and apply \ir{nefcor}(a,c). Similarly, for any rational
$r>\lceil d\sqrt{nl}\rceil$, $H_t$ is nef 
for any rational $r>t>\lceil d\sqrt{nl}\rceil$ by 
\ir{nefcor}(a,c), hence $H_r=(r-t)L'+H_t$ is ample
(since $L'$ meets every curve positively except
for $E_1,\ldots,E_n$, which $H_t$ meets positively).
\qed

\irrnSection{Refinements}{rfnmnts}
The bound $\epsilon(L, n)\ge\varepsilon_{n,l}$
given in \ir{introThm} is limited by the requirement 
in the definition of $\varepsilon_{n,l}$ that $r\le n$.
To get stronger results we need to relax this requirement. 
Our definition of $\varepsilon_{n,l}$ is based on 
\ir{nefcor}, which in turn is based on constructing nef divisors
by blowing up a smooth point of an irreducible curve 
linearly equivalent to a multiple of a very ample
divisor. Considering singular points allows us, in effect,  
to use values of $r$ that can be bigger than $n$.

For example, say $m$ is a positive integer, $p'_1$ is a 
smooth point of an algebraic surface $X$, and 
$X_{p'_1}$ is the blowing up of $X$ at $p'_1$, with $E$ 
being the corresponding exceptional divisor.
If $L$ is very ample on $X$, 
then $tL'-mE=(t-m)L'+m(L'-E)$ is very ample on $X_{p'_1}$
for any $t>m>0$, where $L'$ is the pullback of $L$ to $X_{p'_1}$. 
Thus $|tL'-mE|$ contains an element $C_1$ which 
is reduced and irreducible and is smooth and transverse
to $E$ at some point $p'_2\in E$. Given the morphism
$\pi':Y'\to X$ corresponding to the 
proximity sequence $p'_1,\ldots,p'_n$ with $p'_1$ and $p'_2$
as above and each $p'_i$, $i\le r$, being infinitely near points on proper 
transforms of $C_1$, we find that $[d(\pi'^*L)-mE'_1-E'_2-\cdots-E'_n]$
is the class of an irreducible divisor (in fact,
the proper transform of $C_1$) on $Y'$. 
Define the function $f(d)=\hbox{max}(1,d-1)$;
applying \ir{neflemA}
in the same manner as in \ir{nefcor} we obtain:

\prclm{Corollary}{nefcorRef}{Given $X$, $Y$, $l$, $n$ 
and $L'$ as in \ir{neflem}, let $1\le d$, $1\le m\le f(d)$ 
and $1\le r\le n$ be integers. 
Then we have the following cases.
\item{(a)} If $(r+m-1)^2> nd^2l$, then 
$(r+m-1)L'-dl(E_1+\cdots+E_n)$ is nef.
\item{(b)} If $(r+m-1)^2< nd^2l$, then 
$ndL'-(r+m-1)(E_1+\cdots+E_n)$ is nef.
\item{(c)} If $(r+m-1)^2=nd^2l$, then 
$tL'-(E_1+\cdots+E_n)$ is a nef ${\bf Q}$-divisor
for all rationals $t>\sqrt{n/l}$.}

\noindent If we now define the
sets $$S'_1(n,l)=\Bigl\{\;{r+m-1\over nd} \;\Bigl|\; 1\le r\le n,\; 1\le d, 
\;1\le m\le f(d),\; {r+m-1\over d}\le \sqrt{nl}\;\Bigr\},$$
$$S'_2(n,l)=\Bigl\{\;{dl\over r+m-1}\;\Bigl|\; 1\le r\le n, \;1\le d, 
\;1\le m\le f(d), \;{r+m-1\over d}\ge \sqrt{nl}\;\Bigr\},$$ 
and $S'(n,l)=S'_1(n,l)\cup S'_2(n,l)$, we can take
$\varepsilon'_{n,l}=\hbox{max}(S'(n,l))$. 
Note that since we can rewrite $S'_1$ and $S'_2$ as 
$S'_1(n,l)=\{r/(nd) | 1\le r\le n+f(d)-1, 1\le d, 
r/d\le \sqrt{nl}\}$ and
$S'_2(n,l)=\{dl/r | 1\le r\le n+f(d)-1, 1\le d, 
r/d\ge \sqrt{nl}\}$, this effectively allows us to 
use $r$ bigger than $n$. With essentially the same proof
as for \ir{scgen}, we now have:

\prclm{Theorem}{introThmRef}{Let $l=L^2$, where
$L$ is a very ample divisor on
an algebraic surface $X$. Then $\sqrt{l/n}\ge \epsilon(L, n)$,
and in addition, we have
$\epsilon(L, n)\ge\varepsilon'_{n,l}$
unless $l\le n$ and $nl$ is a square,
in which case $\sqrt{l/n}=\varepsilon'_{n,l}$
and $\epsilon(L, n)\ge\sqrt{l/n}-\varepsilon$
for every positive rational $\varepsilon$.}

\rem{Example}{nplustwoexmpl} This actually is only a minor 
improvement, but it is an improvement.
For example, if $n+2$ is a square, then we can write
$n=s^2+2s-1$ for some $s\le n$. If $s\ge 2$, then 
apply \ir{nefcorRef}(a) with $r=n$, $m=2$ and $d=s+1$
to see that $(r+m-1)L'-d(E_1+\cdots+E_n)$ is nef,
and hence $\epsilon(L, n)\ge\varepsilon'_{n,1}\ge d/(r+m-1)
=(s+1)/(s^2+2s)=\sqrt{1/n}\sqrt{1-{1/(n+1)^2}}$. This is better than what 
we got before (cf. \ir{ptwocor}), and in fact is precisely
the bound obtained in \cite{refBi} for $n=a^2i^2-2i$ for $i=1$.
\vskip\baselineskip

We can get a further effective 
increase in $r$ by considering additional, infinitely near 
singularities. For example, we have:

\prclm{Corollary}{plusonecor}{Say $L$ is a line in $X=\pr2$,
and consider positive integers $d\ge4$, $n\ge5$,
$1\le r'\le n+d-1$. Then, for a blowing up of
$n$ general points of \pr2, $r'L'-d(E_1+\cdots+E_n)$ 
is nef if $r'^2>nd^2$, and $ndL'-r'(E_1+\cdots+E_n)$ 
is nef if $r'^2<nd^2$.}

\noindent{\bf Proof}: If $r'\le n+d-2$, then
we can take $r\le n$ and $m\le d-1$
but still have $r+m-1=r'$, so the result 
follows by \ir{nefcorRef}. Thus we may assume that
$r'=n+d-1$. The idea is to choose a proximity sequence
$p'_1,\ldots,p'_n$ such that
$[dL'-(d-2)E'_1-2E'_2-2E'_3-E'_4-\cdots-E'_n]$ is the class
of an irreducible effective divisor. Given this the result follows
from \ir{neflemA}.

To justify our claim about 
$[dL'-(d-2)E'_1-2E'_2-2E'_3-E'_4-\cdots-E'_n]$,
we first pick $p'_1,\ldots,p'_4$, such that $p'_2$
is on the exceptional divisor of the blow up of $p'_1$,
$p'_3$ is a general point on the exceptional divisor of 
the blow up of $p'_2$ (hence not on the proper transforms 
of the line through $p'_1$ and $p'_2$ nor of the exceptional 
divisor of $p'_1$), and $p'_4$ is a general point on the 
exceptional divisor of the blow up of $p'_3$. The claim is
now that $[dL'-(d-2)E'_1-2E'_2-2E'_3-E'_4]$ is the class
of a reduced irreducible divisor $C_4$ meeting $E_4$ 
transversely. To see this, note that this class corresponds under
a quadratic Cremona transformation centered at $p'_1,p'_2,p'_3$ 
to the class $[(d-2)L''-(d-4)E''_1-E''_4]$, where the $E''_i$
are obtained by blowing up four points with the first three 
as before and the fourth point being a general point on the line
through $p'_1$ and $p'_2$, but not infinitely near to any 
of the first three. But clearly for any $d\ge 4$ there is a reduced
irreducible curve of degree $d-2$ with a $(d-4)$-multiple point
passing simply through some other general point. 
To finish picking our proximity sequence, let $p'_5,\ldots,p'_n$
be the points of the proper transforms 
of $C_4$ infinitely near to $p'_4$.\qed

\rem{Example}{nplusoneexmpl} As an application of the previous result,
let $L$ be a line in \pr2. For $8\le n=s^2+2s$ (thus, $n+1$ is a square), 
we have $\epsilon(L,n)\ge (s^2+3s+1)/(s(s+2)^2)=
\sqrt{1/n}\sqrt{1-{(n-1)/(n(\sqrt{n+1}+1)^2)}}$, and
if $10\le n=s^2+1$ (i.e., $n-1$ is a square), we have 
$\epsilon(L,n)\ge (s+1)/(s^2+s+1))=
\sqrt{1/n}\sqrt{1-{(n-1)/(n+\sqrt{n-1})^2}}$.
To see this, apply \ir{plusonecor}: for $n=s^2+2s$, take
$r=n$, $m=d-2=s$ and $r'=r+m+1$ and note $r'^2<nd^2$, 
and for $n=s^2+1$, take $r=n$, $m=d-2=s-1$ and $r'=r+m+1$
and note $r'^2>nd^2$.
\vskip\baselineskip

We can also obtain additional improvements in
our bounds in special cases, based on the following 
result.

\prclm{Lemma}{adhoc}{Let $d=abc$, where
$a,b,c$ are positive integers with $c<a$
such that $c$ and $a$ are relatively prime
and the characteristic does not divide $c$.
If $r'\ge a^2b^2c$ and $n\ge a^2b^2+(r'-a^2b^2c)$ 
are integers, then, for a blowing up of
$n$ general points of \pr2 with $L\subset\pr2$ 
a line, $r'L'-d(E_1+\cdots+E_n)$ 
is nef if $r'^2>nd^2$, and $ndL'-r'(E_1+\cdots+E_n)$ 
is nef if $r'^2<nd^2$.}

\noindent{\bf Proof}: Let $r=a^2b^2+(r'-a^2b^2c)$.
The idea is to show there is a proximity 
sequence $p_1',\ldots,p_n'$ such that
$dL'-c(E'_1+\cdots+E'_{a^2b^2})-(E'_{a^2b^2+1}+\cdots+E'_r)$
is linearly equivalent to an irreducible divisor,
then apply \ir{neflemA}: if $r'^2>nd^2$, take
$a_0=r'/d$ and $a_1=\cdots=a_n=d$, while
if $r'^2<nd^2$, take $a_0=n$ and $a_1=\cdots=a_n=r'$.

Now we construct our irreducible divisor
$dL'-c(E'_1+\cdots+E'_{a^2b^2})-(E'_{a^2b^2+1}+\cdots+E'_r)$.
We will be very explicit.
Choose homogeneous coordinates $x,y,z$ on \pr2,
let $G=xz^{cb-1}-y^{cb}$, let $F=x^{ab}+z^{(a-c)b}G$.
Note that $G$ and $F$ meet only at $p'_1=[0:0:1]$ (with 
order of contact therefore $ab^2c$), and both
$F$ and $G$ are smooth at $p'_1$. It follows that $F$ and $G$ are reduced
and irreducible, and $p'_1$ is a base point of the pencil
$\langle F^c, G^a\rangle$. This pencil gives a rational map to \pr1.
By successively blowing up points, we can remove the indeterminacies
of this map. It turns out that the points one must blow up
to do so give a proximity sequence $p'_1,\ldots,p'_{a^2b^2}$.
On the blow up the rational map is a morphism, and the 
class of the fiber of the 
morphism corresponding to $F^c$ is just $[cC]=
[abcL'-cE'_1-\cdots-cE'_{a^2b^2}]$, where $C$ is the proper 
transform of the curve defined by $F$. The fiber corresponding
to $G^a$ is $aD+(a-c)N_1+2(a-c)N_2+\cdots+ab^2c(a-c)N_{ab^2c}+
(ab^2c(a-c)-c)N_{ab^2c+1}+
(ab^2c(a-c)-2c)N_{ab^2c+2}+\cdots+cN_{a^2b^2-1}$
(call this divisor $A$ for short)
where $D$ is the proper 
transform of the curve defined by $G$, and $N_i$ is the 
effective divisor whose class is $[E'_i-E'_{i+1}]$. 
Thus $A$ and $cC$ move in a base point free pencil
defining a morphism to \pr1. The divisor $E'_{a^2b^2}$ is a 
multisection of this morphism, since it meets each fiber $c$ times. 
By Bertini's Theorem (see Lemma II.6 of \cite{refars}),
the general member is either reduced and irreducible
or every member is a sum of $c$ elements of $|C|$. 
But the latter would imply
that $A$ is a sum of $c$ members of $|C|$; 
$A$ is connected so $A$ would have to be $c$ times a single
element of $|C|$, which is impossible since $D$ 
is a component of $A$ of multiplicity $a$,
and $a$ and $c$ are relatively prime. Moreover, the trace of
the fibers of the morphism on $E'_{a^2b^2}$ is a linear system
spanned by two points of multiplicity $c$ (since $A$ and $cC$ both meet
$E'_{a^2b^2}$ at single points with multiplicity $c$).
Since the characteristic does not divide $c$, some general fiber $H$ 
is reduced and irreducible and meets $E'_{a^2b^2}$ transversely.
Now take $p'_{a^2b^2+1}$ to be one of these transverse
points of intersection; this uniquely determines the rest of 
the proximity sequence through $p'_r$, with respect to which
$[dL'-c(E'_1+\cdots+E'_{a^2b^2})-(E'_{a^2b^2+1}+\cdots+E'_r)]$
is the class of the proper transform of $H$, which is irreducible.
The rest of the proximity sequence can be chosen
arbitrarily, as long as we don't blow up any more points
of $H$ and keep $E'_i-E'_{i+1}$ irreducible.\qed

\rem{Example}{nminusoneexmpl} Let $n=s^2+j$, where $s$ 
and $j$ are positive integers. If we assume $s$ is not a power of 2
and that the characteristic is not 2, then we may take $c=2$,
$a$ to be any odd prime factor of $s$, $b=s/a$, $d=abc=2s$
and $r'=ca^2b^2+i$, where $i$ is an integer $0\le i\le j$. 
We find $\delta=r'^2-nd^2= a^2b^2c(2i-cj)+i^2$.
This satisfies the hypotheses of \ir{adhoc},
so either $r'L'-d(E_1+\cdots+E_n)$ or $ndL'-r'(E_1+\cdots+E_n)$ 
is nef, depending on the sign of $\delta$.
If we take $i=j$, it follows that 
$r'L'-d(E_1+\cdots+E_n)$ is nef and hence that $\epsilon(L,n)\ge d/r'=
\sqrt{1/n}\sqrt{1-i^2/(2s^2+i)^2}=\sqrt{1/n}\sqrt{1-i^2/(2n-i)^2}$.
When $i=1$ (and hence $n-1$ is a square), this is the bound given 
in \cite{refBi} over ${\bf C}$ (but with no restriction on $s$),
but this remains a very good bound as long as $i$ is not too big.
Similarly, if we take $i=j-1\le2s-1$, then $ndL'-r'(E_1+\cdots+E_n)$ 
is nef and $\epsilon(L,n)\ge r'/(nd)=
\sqrt{1/n}\sqrt{1-(4s^2-i^2)/(4ns^2)}$.
This bound is especially good when $i$ is near $2s$. For example,
if $i=2s-1$ (and hence $n+1$ is a square) we have 
$\epsilon(L,n)\ge\sqrt{1/n}\sqrt{1-(4s-1)/(4ns^2)}>
\sqrt{1/n}\sqrt{1-(4s)/(4ns^2)}=\sqrt{1/n}\sqrt{1-1/(n(\sqrt{n+1}-1))}$.
\vskip\baselineskip

\noindent{\bf Proof} of \ir{introPropB}:
The claims of \ir{introPropB} are proved by 
\ir{nplustwoexmpl}, \ir{nplusoneexmpl} and
\ir{nminusoneexmpl}. \qed

%\vfil\eject
\irrnSection{Applications}{appls}
Our results in \ir{nefblups} have numerous applications
to questions of effectivity, regularity, base point freeness,
ampleness and very ampleness for linear systems on \pr2.

In this section we will always let $L$ be a line in $X=\pr2$ and take 
$\pi:Y\to X$ to be the blow up of $X$ at
$n$ general points $p_1,\ldots,p_n$. Let 
$E_i$, $1\le \ldots\le n$, be the 
corresponding exceptional divisors and let $L'=\pi^*L$.
Given $m>0$, let $F_t=tL'-m(E_1+\cdots+E_n)$.
We can ask:
\item{(a)} What is the least $t$ such that $|F_t|$ is nonempty?
\item{(b)} What is the least $t$ such that $F_t$ is ample?
\item{(c)} What is the least $t$ such that $F_t$ is regular 
(i.e., $h^1(Y,\C O_Y(F_t))=0$)?
\item{(d)} What is the least $t$ such that $|F_t|$ is base point free?
\item{(e)} What is the least $t$ such that $F_t$ is very ample?

\noindent For the rest of this section, 
$F_t$ will be as above.

\irSubsection{Effectivity}{eff}
Here we consider question (a); i.e., what is the least $t$
such that $F_t$ is (linearly equivalent to) 
an effective divisor?

\prclm{Corollary}{effbnds}{If $F_t$ is effective, 
then $t\ge mn\varepsilon_n$.}

\noindent{\bf Proof}: By semicontinuity, $F_t$ remains effective 
under specialization of the points, but 
$N=L'-\varepsilon_n(E_1+\cdots+E_n)$ is nef for some
choice of the points, hence $F_t\cdot N\ge0$. \qed

In terms of simplicity, computability and being 
characteristic free, in addition to its
being a very good bound in an absolute sense,
this bound seems to be the best, overall, now known, at least
for uniform multiplicities. Of course, for 
$D_t=tL'-m_1E_1-\cdots-m_nE_n$ to be effective, it is true that
$t\ge (m_1+\cdots+m_n)\varepsilon_n$, but better bounds
can sometimes be found. For example, if 
$D_t=tL'-2m(E_1+\cdots+E_7)-m(E_8+\cdots+E_{15})$ is effective,
then \ir{nefcorB}(b) (with $d=3$, $j=7$, $r=11$
and $n=15$) gives $t\ge 6m$, whereas the fact that $N$
in the proof of \ir{effbnds} is nef 
gives only $t\ge 22\varepsilon_nm=5.5m$. In some cases of
nonuniform multiplicities, reduction by Cremona transformations
can even give sharp bounds (see \cite{refsurv}). 

Even in the uniform case, there are special cases where
better bounds are known, such as the calculation 
$\epsilon(L, p_1,\ldots,p_{19})\ge 39/170$ in \cite{refBi}
or the examples in \ir{rfnmnts}. 
However, methods which bound effectivity
by testing against nef or ample divisors can at best
say $t>m\sqrt{n}$ if $F_t$ is effective.
Here are some examples of special situations
where better bounds are known:

\item{(a)} Given a positive integer $r$ and $n=4^r$ in characteristic 0,
it follows from \cite{refE} that $F_t$ is effective if and only if
$(t+1)(t+2)>nm(m+1)$ (for $m$ sufficiently large, this is
just $t\ge m\sqrt{n}+(\sqrt{n}-2)/2$). 
\item{(b)} In characteristic 0 when $m\le 12$ and $n>9$, 
\cite{refCM} also proves $F_t$ is effective if and only if
$(t+1)(t+2)>nm(m+1)$.
\item{(c)} The algorithmic bound given
in \cite{refRoe} gives very good bounds, typically better than 
$m\sqrt{n}$, as long as $m$ is not too big compared with 
$n$; however, for $m$ sufficiently large, the
bound in \ir{effbnds} is better (see \cite{refnagprob}). 
\item{(d)} The best overall bounds in characteristic 0 seem to be those
of \cite{refHR}. Although they are asymptotically about the same as those
given here in the sense that they do not seem to lead to better bounds on
Seshadri constants, they do typically give better bounds on effectivity
of $F_t$ for any given $m$. In fact, along with \cite{refE},
\cite{refHR} gives the only bounds currently known which are sharp
in certain cases in which $m$ and $n$ can simultaneously 
(but not independently) be arbitrarily large.

\irSubsection{Ampleness and Regularity}{ampreg}
Here we consider questions (b) and (c); i.e., what is the least $t$
such that $F_t$ is ample, or such that $F_t$ is regular?

\ir{nefblups} already gives a bound for (b):
if $F_c$ is nef, then $F_t$ is ample for all $t>c$. (This is because
$L'$ meets all curves positively except for $E_i$, $i\le n$, but $F_c$
meets each $E_i$ positively.) Consequently we have:

\prclm{Corollary}{ampcor}{If $d>m/\varepsilon_n$,
then $dL'-(E_1+\cdots+E_n)$ is an ample ${\bf Q}$-divisor.}

Now we consider question (c). As was done in \cite{refXb},
duality and the usual vanishing theorems can be used to
convert bounds on nefness or ampleness 
into bounds on regularity. This approach 
gives part (b) of the next result.
(Since \cite{refNtwo} completely solves the regularity problem
for $n\le 9$, we need only consider $n>9$.)

\prclm{Corollary}{reg}{Let $F_t=tL'-m(E_1+\cdots+E_n)$, as usual,
and recall $d^*=\lceil\sqrt{n}\rceil$. Assume $n>9$.
\item{(a)} If $t\ge md^*+\lceil(d^*-3)/2\rceil$,
then $F_t$ is regular.
If $n$ is a square and $m>(d^*-2)/4$, then the converse is true
(i.e., the bound is sharp).
\item{(b)} If $t\ge(m+1)/\varepsilon_n-3$ but $n$ is not a square, 
then $F_t$ is regular.}

\noindent{\bf Proof}: Part (a) follows from Lemma 5.3 of \cite{refHHF}.
For part (b) let $H_{t+3}=F_t-K_Y$; i.e., 
$H_{t+3}=(t+3)L'-(m+1)(E_1+\cdots+E_n)$. Then $H_{t+3}$ is nef and big
(i.e., $H_{t+3}^2>0$) by \ir{nefcor} for 
$t+3\ge(m+1)/\varepsilon_n$, hence 
Ramanujam vanishing (see \cite{reframanujam}, 
or, in positive characteristic, 
Theorem 1.6 of \cite{refterakawa}) implies 
$-H_{t+3}=K_Y-F_t$ is regular,
so by duality $K_Y-(K_Y-F_t)$ is regular. \qed

The bounds given by this corollary seem to be 
the best general bounds now known,
but in special cases better ones are known. For example, 
if $n>9$ is a square but $m$ is not too big, 
the bound in \ir{reg}(a) is known not to be optimal;
in characteristic 0, \cite{refE} gives an optimal bound for all $m$
if $n$ is a power of 4. If $m$ is not too big compared
with $n$, the algorithmic bounds in \cite{refRoeb}, although they are hard
to compute, are often the best available (but for $m$ sufficiently large,
the bounds given here are better; see \cite{refsurv}). 
In \cite{refHR}, bounds are given in characteristic 0
which are better than and sometimes harder to compute
but asymptotically about the same as
those of \ir{reg}(b). The bounds of \cite{refHR} are, however,
sharp for certain values of $m$ and $n$ which can simultaneously 
(but not independently) be
arbitrarily large. Other bounds have also been given.
Those of \cite{refGib}, \cite{refHi} and \cite{refCat} are
on the order of $m\sqrt{2n}$, while those given here are on
the order of $m\sqrt{n}$. Similarly, \ir{reg} is better than
the bound of \cite{refBal} if $m$ is large enough, and
better than \cite{refXb} (Theorem 3) if $n$ is large enough.
See \cite{refsurv} for a discussion and some comparisons.

\irSubsection{Freeness and Very Ampleness}{freeva}
We now consider the last two questions,
what is the least $t$ such that $|F_t|$ is base point free, and
what is the least $t$ such that $F_t$ is very ample?
The results of \ir{ampreg} have an immediate application here.
Indeed, it is well known that $F_t$ is base point free
as long as $F_{t-1}$ is regular, and very ample as long as
$F_{t-1}$ is free and regular. (This follows from the fact that the ideal
$I_Z$ of the fat point subscheme $Z=m(p_1+\cdots+p_n)$
is generated in degrees $t\le\sigma$ \cite{refDGM}, where $\sigma$
can be defined as one more than the least $t$ such that 
$F_t$ is regular.) Thus, if $F_t$ is regular for all
$t\ge N$ for some $N$, then $|F_t|$ has no base points 
for $t\ge N+1$ and is very ample for $t\ge N+2$.

In certain cases, one can do better. 
For example, \cite{refHHF} shows that when $n>9$ is an even square and
$m>(\sqrt{n}-2)/4$, then $F_t$ is both regular for all 
$t\ge m\sqrt{n}+(\sqrt{n}-2)/2$ and that $I_Z$ is generated in degrees 
at most $m\sqrt{n}+(\sqrt{n}-2)/2$. Thus $F_t$ 
is also base point free for all 
$t\ge m\sqrt{n}+(\sqrt{n}-2)/2$, and very ample for all 
$t\ge m\sqrt{n}+\sqrt{n}/2$. For additional (but characteristic 0) 
examples, when $n$ is not a square but both $n$ and $m$ can be  
large, see \cite{refHR}.

\References

\bibitem{\refBal} E. Ballico. {\it Curves of minimal degree
with prescribed singularities}, Illinois J.\ Math.\ 45 (1999), 672--676.

\bibitem{\refBau} T. Bauer. {\it Seshadri constants on algebraic surfaces}, 
Math.\ Ann.\ 313 (1999), 547--583.

\bibitem{\refBi} P. Biran. {\it Constructing new ample 
divisors out of old ones}, Duke Math. J. 98 (1999), no. 1, 113--135.

\bibitem{\refCat} M.\ V. Catalisano. {\it Linear Systems of Plane Curves 
through Fixed ``Fat'' Points of \pr2}, 
J.\ Alg.\  142 (1991), 81-100.

\bibitem{\refCM} C. Ciliberto and R. Miranda. {\it Linear systems
of plane curves with base points of equal multiplicity}, 
Trans. Amer. Math. Soc. 352 (2000), 4037--4050.

\bibitem{\refDGM} E.\ D. Davis, A.\ V. Geramita, and P.
Maroscia. {\it Perfect
Homogeneous Ideals: Dubreil's Theorems Revisited},
Bull.\ Sc.\ math., $2^e$ s\'erie, 108 (1984), 143--185.

\bibitem{\refdem} J.-P. Demailly. {\it Singular Hermitian metrics on
positive line bundles}, Complex Algebraic Varieties, Proceedings,
1990, LNM 1507 (1992), 87--104.

\bibitem{\refEL} L. Ein and R. Lazarsfeld.
{\it Seshadri constants on smooth surfaces},
Journ\'ees de G\'eom\'etrie Alg\'ebrique d'Orsay, Ast\'erisque 282
(1993), 177-186.

\bibitem{\refE} L. Evain.
{\it La fonction de Hilbert de la r\'eunion
de $4^h$ gros points g\'en\'eriques
de \pr2 de m\^eme multiplicit\'e},
J. Alg. Geom. 8 (1999), 787--796.

\bibitem{\refGib} A. Gimigliano.
{\it Regularity of Linear Systems of Plane Curves},
J. Alg. 124 (1989), 447--460.

\bibitem{\refHHF} B. Harbourne, S. Holay and S. Fitchett.
{\it Resolutions of Ideals of Quasiuniform Fat Point Subschemes of \pr2},
preprint (2000).

\bibitem{\refars} B. Harbourne. {\it  Anticanonical Rational surfaces},
Trans. Amer. Math. Soc. 349 (1997), 1191--1208.

\bibitem{\refnagprob} \manyby. {\it  On Nagata's Conjecture},
J. Alg. 236 (2001), 692--702.

\bibitem{\refsurv} \manyby. {\it  Problems and Progress:
A survey on fat points in \pr2},
preprint (2000).

\bibitem{\refHR} B. Harbourne and J. Ro\'e. 
{\it Linear systems with multiple base points in \pr2}, preprint (2000).

\bibitem{\refHi} A. Hirschowitz.
{\it Une conjecture pour la cohomologie 
des diviseurs sur les surfaces rationelles g\'en\'eriques},
Journ.\ Reine Angew.\ Math. 397
(1989), 208--213.

\bibitem{\refKu} M. K\"uchle.
{\it Multiple point Seshadri constants and the dimension
of adjoint linear series},
Ann.\ Inst.\ Fourier, Grenoble, 46 (1996), 63--71.

\bibitem{\refKub} \manyby.
{\it Ample line bundles on blown up surfaces},
Math.\ Ann.\ 304 (1996), 151--155.

\bibitem{\refNone} M. Nagata. {\it On the 14-th problem of Hilbert}, 
Amer.\ J.\ Math.\ 81 (1959), 766--772.

\bibitem{\refNtwo} \manyby. {\it On rational surfaces, II}, 
Mem.\ Coll.\ Sci.\ 
Univ.\ Kyoto, Ser.\ A Math.\ 33 (1960), 271--293.

\bibitem{\reframanujam} C. P. Ramanujam. {\it Supplement
to the article ``Remarks on the Kodaira
vanishing theorem''}, J. Indian Math. Soc. {\bf 38} (1974), 121--124.

\bibitem{\refRoe} J. Ro\'e. {\it On the existence of plane curves
with imposed multiple points}, J. Pure Appl. Alg. 156(2001), 115--126.

\bibitem{\refRoeb} \manyby. {\it Linear systems of plane curves with
imposed multiple points}, preprint (2000).

\bibitem{\refS} T. Szemberg. {\it Global and local positivity of 
line bundles}, Habilitation, 2001.

\bibitem{\refST} T. Szemberg and H. Tutaj-Gasi\'nska. {\it General
blow ups of the projective plane}, to appear, Proc. Amer. Math. Soc.

\bibitem{\refterakawa} H. Terakawa. {\it The $d$-very ampleness
on a projective surface in characteristic $p$}, 
Pac. J. Math. 187 (1999), 187--199.

\bibitem{\refXa} G. Xu. {\it Divisors on the blow up of the projective 
plane}, Man.\ Math.\ 86 (1995), 195--197.

\bibitem{\refXb} \manyby. {\it Ample line bundles on smooth surfaces}, 
Jour.\ Reine Ang.\ Math.\ 469 (1995), 199--209.

\bye